\documentclass[a4paper, 10pt]{article}
\usepackage[english]{babel}
\usepackage[utf8]{inputenc}
\usepackage[T1]{fontenc}
\usepackage[leqno]{amsmath}
\usepackage{amssymb,amsfonts,mathrsfs,exscale,colonequals,wasysym}
\usepackage[linkcolor=blue]{hyperref}
\usepackage{graphicx}
\usepackage{mathtools}
\usepackage{array}
\usepackage{longtable}
\usepackage{geometry}
\usepackage{pdfpages}
\usepackage{graphicx}
\usepackage{tikz-cd}

\def\bigtimes{\mathop{\raisebox{-2pt}{\huge $\times$\kern-1pt}}}

\def\bigboxplus{\mathop{\raisebox{-2pt}{\scalebox{1.6}{$\kern1pt\boxplus\kern1pt$}}}}

\def\Spec{\mathop{\rm Spec}\nolimits}

\def\rank{\mathop{\rm rank}\nolimits}

\def\charact{\mathop{\rm char}\nolimits}

\def\id{{\rm id}}

\let\phi\varphi
\let\theta\vartheta
\let\epsilon\varepsilon
\let\setminus\smallsetminus

\let\leq\leqslant
\let\geq\geqslant

\usepackage{centernot}

\newcommand{\BA}{{\mathbb{A}}}

\newcommand{\BN}{{\mathbb{N}}}

\newcommand{\BR}{{\mathbb{R}}}

\newcommand{\BZ}{{\mathbb{Z}}}

\newcommand{\Fm}{{\mathfrak{m}}}

\newcommand{\CI}{{\cal I}}

\newcommand{\CQ}{{\cal Q}}

\usepackage{amsthm}

\makeatletter
\renewcommand\thm@space@setup{%
  \thm@preskip=8pt plus 2pt minus 2pt%
  \thm@postskip=2pt plus 1pt minus 1pt%
}
\makeatother

\newtheorem{Thm}{Theorem}[subsection]
\newtheorem{Lem}[Thm]{Lemma}
\newtheorem{Prop}[Thm]{Proposition}
\newtheorem{Cor}[Thm]{Corollary}

\theoremstyle{definition}
\newtheorem{Ex}[Thm]{Example}
\newtheorem*{Thm*}{Theorem}

\newtheorem{Def}[Thm]{Definition}

\theoremstyle{remark}

\newtheorem*{Ques}{Open Question}
\newtheorem*{Conj}{Conjecture}

\def\qed{{\hskip0pt\unskip\unskip\nobreak\hfil\penalty50
          \hskip1em\hbox{}\nobreak\hfil
           {$\square$}
          \parfillskip=0pt\finalhyphendemerits=0
          \par}}

\newenvironment{myproof}
               {\noindent{\textit{Proof}.}\;}
               {\qed}

\newenvironment{skproof}
               {\noindent{\textit{Sketch of a proof}.}\; }
               {\qed}

\numberwithin{equation}{section}

\usepackage{enumitem} % control layout of enumerate, itemize
\setlist{nolistsep}
{\begin{enumerate}[label=(\alph*),topsep=0pt,partopsep=10pt,itemsep=1ex]}%
{\end{enumerate}}
{\begin{enumerate}[label=(\roman*),topsep=0pt,partopsep=5pt,itemsep=0.5ex]}{\end{enumerate}}
\newenvironment{myitemize}{\begin{itemize}[label=-]}{\end{itemize}}

\usepackage{appendix}

\newcommand{\ideal}{\trianglelefteq}

{\makeatletter
\g@addto@macro\normalsize{%
  \setlength\abovedisplayskip{10pt}
  \setlength\belowdisplayskip{8pt}
  \setlength\abovedisplayshortskip{-5pt}
  \setlength\belowdisplayshortskip{5pt}
}
\makeatother}

\title{A Complete Classification of Ideal Chomp Games on Low-Rank Algebras
}
    \author{Leopold Karl}
    \date{August 2025.}

\begin{document} 

\maketitle
\thispagestyle{empty}

%%%%%%%%%%%%%%%%%%%%%%%%%%%%%%%%%%%%%%%%%%%%%%%%%%%%%%%%%%%%%%%%%%%%%%

\section*{Abstract}
%\addcontentsline{toc}{section}{Abstract}%

%Abstract

We completely classify winning strategies in the Ideal Chomp Game played on $\bar{K}$-algebras $R$ of rank at most 6. In this two-player combinatorial game, players 
alternately add generators to build an ideal inside a given ring $R$, with the player who builds an ideal equal to the entire ring losing. We prove that player A has a winning strategy on all $\bar{K}$-algebras $R$ up to rank 6 except for five specific cases: $\bar{K}$ itself, $\bar{K}[x,y]/(x,y)^2$, and three other local algebras. Our methods combine game-theoretic analysis with the structure theory of Artinian rings and computational verification. We also discuss a classical result of Henson on winning strategies in the Ideal Chomp Game, as well as ideas and open questions about the Ideal Chomp Game on higher-dimensional $\bar{K}$-algebras.

\tableofcontents
\pagebreak

%%%%%%%%%%%%%%%%%%%%%%%%%%%%%%%%%%%%%%%%%%%%%%%%%%%%%%%%%%%%%
% Chapters

\section{Introduction}

% Motivation

\begin{figure}[h!]
% \label{ChompGame}
    \centering
    \includegraphics[width=0.9\textwidth]{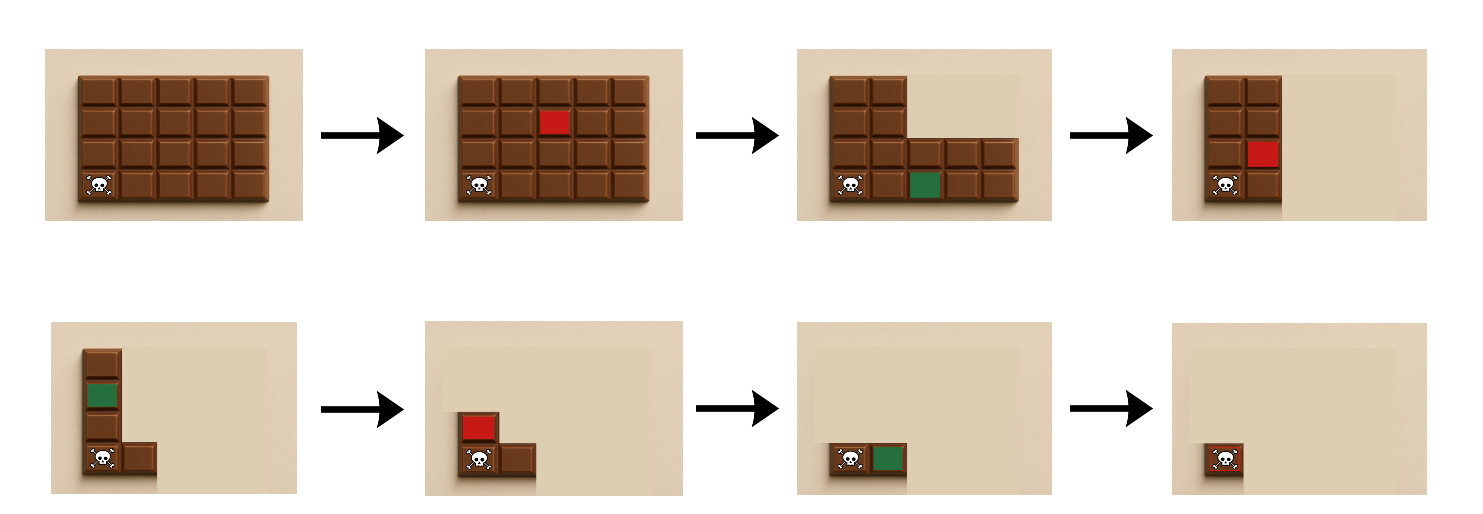}  % Adjust path and width as needed
    \caption{A game of Chomp where player A (red) loses.}
    \label{fig:chomp}
\end{figure}

In this paper, we study the \emph{Ideal Chomp Game}, a two-player game whose ``playing field'' is a given fixed ring $R$. The origins of the Ideal Chomp Game lie in the (classical) Chomp Game formulated by David Gale \cite{gale} in 1974. In the classical Chomp Game two players, Alice (A) and Bob (B) chomp pieces of a rectangular chocolate bar. More precisely, they choose a piece and chomp all pieces above and to the right of it (see Figure \ref{fig:chomp}). The player who eats the last piece of chocolate, which according to the rules of the game is always the lowest and leftmost piece, loses the game.

The Ideal Chomp Game is in spirit similar to the Chomp Game, but played on the algebraic structure of a ring, which allows for a geometric interpretation. It was first studied by C.W. Henson in his paper \cite{henson} in 1970 and is played according to the following set of rules: First, a ring $R$ is fixed on which the two players $A$ and $B$ decide to play the Ideal Chomp Game. The players then start with the zero ideal $I_0 = (0)$ and in each turn increase the ideal $I_{n-1}$ strictly by adding a generator to obtain $I_{n}$. The player who increases the ideal to the entire ring loses. On the geometric side, the move of adding the element $a$ to the set of generators corresponds to the intersection of the current subscheme $$S_{n-1} = \Spec R/I_{n-1} \hookrightarrow \Spec R$$ corresponding to the ideal $I_{n-1}$ with the subscheme corresponding to the ideal $(a)$ . The player who turns $S_n$ into the empty subscheme, i.e. takes the last point away, loses the game. Notice the similarity to the classical Chomp Game.

Indeed, the Classical Chomp Game on a rectangle of size $a \times b$ is essentially the same as the Ideal Chomp Game on the ring $K[x,y]/(x^a, y^b)$ if we additionally restrict the possible moves of the players from polynomials to monomials.

\begin{Ex}
    On $\BZ$ a possible game could look like the following:
    First, Alice plays $f_1 = 100$. Next, Bob plays $f_2 = 2 \in \BZ \setminus (100)$. Now whatever element $f_3 \in \BZ \setminus (2)$ Alice chooses, she will lose since $(f_1,f_2,f_3)$ will generate $\BZ$. Clearly, Alice could have won the Ideal Chomp Game on $\BZ$ by starting with a prime.
\end{Ex}

The main result of our paper is the following complete classification of the Ideal Chomp Game on $\bar{K}$-algebras $R$ up to rank $6$.

%%%
\iffalse
\begin{Ex}
    On the ring $R = K[x,y]/(x,y)^2$ the game could play out as follows: First, Alice chooses $f_1 = x + y \in R \setminus \{0\}$. Then Bob chooses $f_2 = 2x + y \in R \setminus (x+y)$, so $S_2 = (x,y)$. No matter what Alice plays next, she loses.
\end{Ex}
\fi
%%%

\enlargethispage{3mm}

\begin{Thm}
    [Classification of the Ideal Chomp Game on $\bar{K}$-algebras up to rank 6]
    \label{ClassificationThm1}
    Alice has a winning strategy in the Ideal Chomp Game on any $\bar{K}$-algebra $R$ up to rank 6, except if $R$ is isomorphic to \vspace{-3mm}
    \begin{align*}
        R_1 & = \bar{K},\\
        R_4 & = \bar{K}[x,y]/(x,y)^2,\\
        R_{12} & =  \bar{K}[x,y]/(xy,x^3,y^3),\\
        R_{13} & = \bar{K}[x,y]/(x^2,xy^2, y^3),\\
        R_{17} & = \bar{K}[x,y,z,w]/(x,y,z,w)^2.
    \end{align*}
\end{Thm}

\pagebreak

\section{Preliminaries}

% Game Description
\subsection{Game Instructions}
We proceed with a more formal description of the Ideal Chomp Game.

\begin{Def}
    [Ideal Chomp Game] Let $R$ be a ring and $I_0 = (0)$ the zero ideal. Players A and B alternately take turns, starting with player A. In each turn, the current player chooses an element $a_n \in R \setminus I_{n-1}$ and builds the ideal $I_n = I_{n-1} + (a_n)$. The player who increases the ideal to the entire ring, that is, who sets $I_n = R$, loses the game. 
\end{Def}

More explicitly, the $n$-th ideal is given by $I_n = (0, a_1, a_2, \dots, a_n)$. In addition, a valid state of an Ideal Chomp Game is completely described by the tuple $(R, a_1, \dots, a_N)$ where $R$ is the ring on which the game is played and $a_1, \dots, a_N$ are the elements the players choose in the subsequent turns, i.e. satisfy $a_i \in R \setminus (0,a_1,\dots,a_{i-1})$.

\begin{Ex}
    The Ideal Chomp Game $(\BR[x,y], x+y, 2x + y, 2)$, is a valid game since
$$x+y \in \BR[x,y],$$ 
$$2x+y \in \BR[x,y]\setminus (x+y) \text{ and}$$ 
$$2 \in \BR[x,y] \setminus (x+y, 2x+y).$$
The game is won by the second player since the first player increased the ideal to the entire ring by adding $2$: $(x+y,2x + y,2) = (x,y,2) = \BR[x,y].$
\end{Ex}

% Reformulation
\subsection{A Reformulation of the Game}
The Ideal Chomp Game can also be phrased in the following way.

\begin{Def}[Quotient Chomp Game]

    Let $R_0 = R$ be a ring. Players A and B alternately take turns, starting with player A. In each turn, the current player chooses an element $a_n \in R_{n-1} \setminus \{0\}$ and builds the ring $R_{n} = R_{n-1} / (a_n)$. The player who reduces the given ring to the zero ring, that is, who sets $R_{n} = 0$, loses the game.
\end{Def}

More explicitly, the $n$-th ring is given by $$R_n = \left(((R/a_1)\dots)/a_n\right) \cong R/(a_1,\dots,a_{n}).$$ In addition, the state of a game is completely described by the tuple $(R, a_1,\dots, a_n)$ where $a_i \in R_{i-1}$.

Due to the following proposition, the Quotient Chomp Game is, indeed, only a reformulation of the Ideal Chomp Game.

\begin{Prop}
    The Ideal Chomp Game on the ring $R$ and the Quotient Chomp Game starting with the ring $R_0 = R$ are equivalent formulations of the same game.
\end{Prop}

\begin{myproof}
    %Consider the Ideal Chomp Game on the ring $R$ and the Quotient Chomp Game starting on the ring $R_0 = R$. We show that the game states of these two games stand in bijection and that player X, with $X \in \{A,B\}$, has a winning strategy in the Ideal Chomp Game if and only if he has a winning strategy in the Quotient Chomp Game.
    %For the Ideal Chomp Game we call a tuple $(R,a_1,\dots,a_N)$ where the first entry specifies the ring the game is played on and the $a_i$'s are a sequence of allowed moves such that the game ends with $a_N$, a complete set of moves.
    Let $\CI_R = \{(R,a_1,\dots, a_N) \mid \forall n \in [N]: a_n \notin (a_1,\dots, a_{n-1})\}$ be the set of all complete sets of moves in the Ideal Chomp Game on the ring $R$. Since two complete sets of moves reflect the same game if and only if all the ideals $I_n = (a_1,\dots,a_n)$ equal one another, we introduce the equivalence relation $$(R,a_1,\dots, a_N) \sim (R, b_1, \dots, b_N) \iff \forall n \in [N]: (a_1, \dots, a_n) = (b_1,\dots, b_n).$$
    
    Then the set of all valid Ideal Chomp Games is $\CI^{\sim}_R = \CI_R/\text{\hspace*{-1.5mm}}\sim$. 
    
    Similarly, let $\CQ_R = \{(R,a_1,\dots, a_N) \mid \forall n \in [N]: a_n \in R/(a_1,\dots, a_{n-1})\setminus\{0\}\}$ be the set of all complete sets of moves in the Quotient Chomp Game starting on the ring $R_0 = R$ and let $\CQ_R^{\sim} = \CQ_R/\text{\hspace*{-1.5mm}}\sim$ be the set of all valid Quotient Chomp Games where the equivalence relation is given as
    \[(R,a_1,\dots, a_N) \sim (R,b_1,\dots,b_N) \iff \forall n \in [N]\ \exists \phi_n:
    \begin{tikzcd}
    R \arrow[r]\arrow[dr] & R_{a,n} \arrow[d, "\phi_n"] \\
     & R_{b,n}
    \end{tikzcd}\
    %R_{a,n} = R_{b,n}
    \text{commutes.}\]
    where $R_{i,n}$ for $i \in \{a,b\}$ is inductively defined by $R_{i,0} = R$ and $R_{i,n} = R_{i,n-1}/(i_n)$.

    We claim that the maps $$\phi: \CI_R^{\sim} \to \CQ_R^{\sim}, [(R, a_1, \dots, a_N)] \mapsto [(R, [a_1], \dots, [a_N])] \hspace*{5mm} \text{and}$$ $$\psi: \CQ_R^{\sim} \to \CI_R^{\sim}, [(R, b_1, \dots, b_N)] \mapsto [(R, \psi(b_1), \dots, \psi(b_N))],$$
    where $[a_i]$ is the equivalence class of $a_i$ in $R_{a,i-1}$ and $\psi(b_i) \in R$ is any element with $[\psi(b_i)] = b_i$ in $R_{i-1}$, are well-defined and two-sided inverses to one another.

    For a valid game $[(R, a_1, \dots, a_N)] \in \CI_R^{\sim}$, the game $\phi([(R,a_1,\dots,a_N)])$ is a well-defined element of $\CQ_R^{\sim}$ as the rings $R_{a,n}$ only depend on the ideals $I_n$, i.e. only on the game's equivalence class. For a valid game $[(R, b_1, \dots, b_N)] \in \CQ_R^{\sim}$ we have that for any $i$ and any choice of representative $\psi(b_i)$ it holds that $$\psi(b_i) \in R \setminus (0,\psi(b_1), \dots, \psi(b_{i-1}))$$ since $b_i \neq 0$ in $R_{i-1}$. In addition, the game $[(R, \psi(b_1),\dots, \psi(b_N))]$ does not depend on the choice of representative $\psi(b_i)$ since for the cosets of any two representatives $x, y \in b_i$ we have $[x] = [y]$ in $R/(0,\psi(b_1), \dots, \psi(b_{i-1}))$, so the ideals $(0,\psi(b_1), \dots, \psi(b_{i-1}), x)$ and $(0, \psi(b_1), \dots, \psi(b_{i-1}), y)$ lie in the same equivalence class of $\CI_R^{\sim}$, that is, they are representatives of the same element in $\CI_R^{\sim}$. Furthermore, there exists an element $x_{b_i} \in R$ that satisfies $[x_{b_i}] = b_i$ in $R_{b,i-1}$ since $\pi_{i-1}^{-1}(b_i)$ is non-empty where $$\pi_{i-1}: R \to R_{b,i-1} = \left((R/(b_0))/\dots\right)/(b_{i-1})$$ is the projection map.

    To show that the maps $\phi$ and $\psi$ are two-sided inverses of one another, notice that if we set $I_{i-1} = (0,a_1, \dots, a_{i-1})$, then $a_i$ is a representative of the coset $[a_i] = a_i I_{i-1}$. Since $\psi$ is independent of the choice of representative by the above paragraph, we conclude that $\psi \circ \phi = \id_{\CI_R^{\sim}}$. The direction $\phi \circ \psi = \id_{\CQ_R^{\sim}}$ is immediate as the coset of a representative of a coset is again the same coset.

    Finally, a valid Ideal Chomp Game is won by the first player if and only if the corresponding valid Quotient Chomp Game is won by the first player since the number of turns is the same and both games are won by the first player if and only if the game consists of an even number of turns. Similarly, the second player wins the game if and only if the game consists of an odd number of turns. Hence, these games are just reformulations of one another.
\end{myproof}

% First Winning Conditions
\subsection{First Winning Conditions}

In the following, we list some immediate winning conditions.

\begin{Prop}
\label{fieldwin}
    If $R$ is a field, then player $B$ wins.
\end{Prop}

\begin{myproof}
    Already in the first round, player $A$ has to increase the zero ideal to the entire ring $R$ since the only proper ideal of a field is the zero ideal. Player $B$ wins although they did not even make a move.
\end{myproof}

\begin{Prop}
\label{principalmaxidealwin}
    If $R$ is a ring with a principal maximal ideal and not a field, then player $A$ wins. In particular, player $A$ wins on any principal ideal domain that is not a field.
\end{Prop}

\begin{myproof}
    Let $(a) \ideal R$ be a principal maximal ideal. Since $R$ is not a field, $(a) \neq (0)$, the move $I_1 = (a)$ is a regular move for player $A$. Now, player $B$ has to increase the maximal ideal to the entire ring. Hence, player $A$ wins the game.
\end{myproof}

Notice that at least for a Notherian ring $R$ the question ``Which player has a winning strategy in the Ideal Chomp Game on $R$?'' is well-posed due to the following corollary.

\begin{Cor}\label{ExistenceOfWinningStrategy}
    Let $R$ be a Noetherian ring. Then one of the players has a winning strategy in the Ideal Chomp Game.
\end{Cor}

\begin{myproof}
    By Zermelo's Theorem for deterministic, finite, two-player, zero-sum games of perfect information (see \cite{zermelo}), there exist strategies such that the game ends in a draw or one of the players has a winning strategy. Since the Noetherian condition implies the finiteness of the Ideal Chomp Game and there is no mechanism for a game ending in a draw in the Ideal Chomp Game, we conclude that one of the players has to have a winning strategy in the Ideal Chomp Game on $R$.
\end{myproof}

Finally, we can easily determine the winner of the Ideal Chomp game on any Cartesian product of rings.

\begin{Prop}
\label{CartesianProdCor}
    Let $R = \prod_{i \in I} R_i$ be a non-trivial Cartesian product of Noetherian rings, i.e. for all $i \in I$ we have $R_i \neq 0$ and $|I| \geq 2$. Then player A has a winning strategy in the Ideal Chomp Game on $R$.
\end{Prop}

\begin{myproof}
    If there exists an index $i_A \in I$ such that player A has a winning strategy on $R_{i_A}$, then player A has a winning strategy on $R$ by reducing the game on R to the game on $R_{i_A}$ by playing $a = (a_i)_{i \in I}$ with $$a_i = \begin{cases}
        a_W & \text{if } i = i_A\\
        1 & \text{otherwise}
    \end{cases}$$ in her first move where $a_W$ is the element player A's winning strategy on $R_{i_A}$ starts with.\newline Otherwise, the second player has a winning strategy on every $R_i$ by Corollary \ref{ExistenceOfWinningStrategy}. In particular, player B has a winning strategy in the Ideal Chomp Game on $R_1$. Player A can steal this winning strategy of player B on $R_1$ by playing $a = (a_i)_{i \in I}$ with $$a_i = \begin{cases}
        0 & \text{if } i = 1\\
        1 & \text{otherwise}
    \end{cases}$$
    which leaves player B with the ring $R_1$ on which the second player (now player A) has a winning strategy.
\end{myproof}

% Poonen's classification
\subsection{Poonen's Classification of local $\overline{K}$-Algebras up to Rank 6}

From the perspective of the Quotient Ring Game it is clear that the dimension of the ring $R$ is monotonically decreasing throughout the game play. This suggests the approach where one tries to classify rings of increasing dimension one after the other by finding second-player-win rings $S$ and noticing that any ring $R$ that player A can reduce to $S$ by factoring through a principal ideal is a first-player-win ring. However, even the case of zero-dimensional rings is quite tricky already. The most basic examples are $\bar{K}$-algebras R which are finite dimensional as $\bar{K}$-vector spaces and, hence, also Artinian rings. By the Structure Theorem for Artinian rings (see Theorem 8.7 in \cite{atiyah-macdonald}), any such $\bar{K}$-algebra $R$ is a finite product of local $\bar{K}$-algebras. By Proposition \ref{CartesianProdCor} immediately classify all $\bar{K}$-algebras with more than one factor in this factorisation. Hence, we are left with classifying the Ideal Chomp Game on local $\bar{K}$-algebras $R$.

In this case, Poonen classified all local $K$-algebras up to rank $6$ in his paper \cite{poonen}. We summarize his results in the following Table \ref{poonentable}.
Here $n$ is the rank of the $K$-algebra, the vector $\vec{d}$ is given as $\vec{d} = (d_i)_{i > 0} = (\dim(\Fm^{i}/\Fm^{i+1}))_{i > 0}$ where $\Fm$ is the unique maximal ideal of the local $K$-algebra.  An asterisk ``$*$'' means that this case is only a separate isomorphism class when $\charact(K) = 2$ and a double asterisk ``$**$'' means that this case is only a separate isomorphism class if $\charact(K) = 3$, in other characteristica these cases reduce to the case immediately above it. The column ``Win'' summarises the results we will obtain throughout the next Section \ref{ClassificationIdealChomp} by showing which player has a winning strategy on the respective $K$-algebra.\\

%\vspace{3cm}

%%%%%%%%% Table
%\enlargethispage{1cm}

\renewcommand{\arraystretch}{1.2}
\setlength\LTleft{-1cm}
\setlength\LTright{0pt plus 1fill}
\thispagestyle{empty}
\enlargethispage{10mm}

\begin{longtable}{|c|c|c|c|c|}
        \hline
        \( n \) & \( \vec{d} \) & Name & \textbf{Local $K$-algebra} & \textbf{Win}\\
        \hline
        1 & (0) & \(R_1\) & \( K \) & B\\
        \hline
        2 & 1 & \(R_2\) & \( K[x] / (x^2) \) & A\\
        \hline
        3 & 1,1 & \(R_3\) & \( K[x] / (x^3) \) & A\\
          & 2 & \(R_4\) & \( K[x, y] / (x, y)^2 \) & B\\
        \hline
        4 & 1,1,1 & \(R_5\) & \( K[x] / (x^4) \) & A\\
          & 2,1 & \(R_6\) & \( K[x, y] / (x^2, xy, y^3) \) & A\\
          & 2,1 & \(R_7\) & \( K[x, y] / (x^2, y^2) \) & A\\
          & *2,1 & \(R_{7,*}\) & \( K[x, y] / (x^2 + y^2, xy) \) & A \\
          & 3 & \(R_8\) & \( K[x, y, z] / (x, y, z)^2 \) & A\\
        \hline
        5 & 1,1,1,1 & \(R_9\) & \( K[x] / (x^5) \) & A\\
          & 2,1,1 & \(R_{10}\) & \( K[x, y] / (x^2, xy, y^4) \) & A\\
          & 2,1,1 & \(R_{11}\) & \( K[x,y] / (x^2 + y^3, xy) \) & A\\
          & 2,2 & \(R_{12}\) & \( K[x, y] / (xy, x^3, y^3) \) & B\\
          & 2,2 & \(R_{13}\) & \( K[x,y] / (x^2, xy^2, y^3) \) & B\\
          & 3,1 & \(R_{14}\) & \( K[x, y, z] / (x^2, y^2, xy, xz, yz, z^3) \) & A\\
          & 3,1 & \(R_{15}\) & \( K[x, y, z] / (x^2, y^2, z^2, xy, xz) \) & A\\
          & *3,1 & \(R_{15,*}\) & \( K[x, y, z] / (x^2, xy, xz, yz, y^2 + z^2) \) & A\\
          & 3,1 & \(R_{16}\) & \( K[x, y, z] / (xy, xz, yz, x^2 + y^2, x^2 + z^2) \) & A\\
          & 4 & \(R_{17}\) & \( K[x, y, z, w] / (x, y, z, w)^2 \) & B\\
        \hline
        6 & 1,1,1,1,1 & \(R_{18}\) & \( K[x] / (x^6) \) & A\\
          & 2,1,1,1 & \(R_{19}\) & \( K[x, y] / (x^2, xy, y^5) \) & A\\
          & 2,1,1,1 & \(R_{20}\) & \( K[x,y] / (x^2+y^4, xy) \) & A\\
          & 2,2,1 & \(R_{21}\) & \( K[x, y] / (xy, x^3, y^4) \) & A\\
          & 2,2,1 & \(R_{22}\) & \( K[x,y] / (xy, x^3 + y^3) \) & A\\
          & 2,2,1 & \(R_{23}\) & \( K[x,y] / (x^2, xy^2, y^4) \) & A\\
          & 2,2,1 & \(R_{24}\) & \( K[x,y] / (x^2+y^3, xy^2, y^4) \) & A\\
          & 2,2,1 & \(R_{25}\) & \( K[x,y] / (x^2, y^3) \) & A\\
          & *2,2,1 & \(R_{25,*}\) & \( K[x,y] / (x^2 + xy^2, y^3) \) & A\\
          & **2,2,1 & \(R_{25,**}\) & \( K[x,y] / (x^2, xy^2 + y^3) \) & A\\
          & 2,3 & \(R_{26}\) & \( K[x,y] / (x,y)^3 \) & A\\
          & 3,1,1 & \(R_{27}\) & \( K[x, y, z] / (x^2, xy, y^2, xz, yz, z^4) \) & A\\
          & 3,1,1 & \(R_{28}\) & \( K[x,y,z] / (x^2,xy,y^2+z^3, xz, yz, z^4) \) & A\\
          & 3,1,1 & \(R_{29}\) & \( K[x,y,z] / (x^2, xy + z^3, y^2, xz, yz, z^4) \) & A\\
          & * 3,1,1 & \(R_{29,*}\) & \( K[x,y,z] / (x^2 + z^3, xy, y^2 + z^3, xz, yz, z^4) \) & A\\
          & 3,2 & \(R_{30}\) & \( K[x, y, z] / (xy, yz, z^2, y^2 - xz, x^3) \) & A\\
          & 3,2 & \(R_{31}\) & \( K[x,y,z] / (xy, z^2, xz - yz, x^2 + y^2 - xz) \) & A\\
          & *3,2 & \(R_{31,*}\) & \( K[x, y, z] / (x^2, z^2, y^2 - xz, yz) \) & A\\
          & 3,2 & \(R_{32}\) & \( K[x,y,z] / (x^2, xy, xz, y^2, yz^2, z^3) \) & A\\
          & 3,2 & \(R_{33}\) & \( K[x,y,z] / (x^2, xy, xz, yz, y^3, z^3) \) & A\\
          & 3,2 & \(R_{34}\) & \( K[x,y,z] / (xy, xz, y^2, z^2, x^3) \) & A\\
          & *3,2 & \(R_{34,*}\) & \( K[x,y,z] / (xy, xz, yz, y^2- z^2, x^3) \) & A\\
          & 3,2 & \(R_{35}\) & \( K[x,y,z] / (xy, xz, yz, x^2 + y^2 - z^2) \) & A\\
          & 3,2 & \(R_{36}\) & \( K[x,y,z] / (x^2, xy, yz, y^2 - z^2) \) & A\\
          & *3,2 & \(R_{36,*}\) & \( K[x,y,z] / (x^2, xy, yz, xz + y^2 - z^2) \) & A\\   
          & 3,2 & \(R_{37}\) & \( K[x,y,z] / (x^2, xy,y^2, z^2) \) & A\\
          & *3,2 & \(R_{37,*}\) & \( K[x,y,z] / (x^2, xy, y^2, z^2 - xz) \) & A\\
          & 4,1 & \(R_{38}\) & \( K[x, y, z, w] / (x^2, y^2, z^2, xy, xz, xw, yz, yw, zw, w^3) \) & A\\
          & 4,1 & \(R_{39}\) & \( K[x,y,z,w] / (x^2, y^2, z^2, w^2, xy, xz, xw, yz, yw) \) & A\\
          & *4,1 & \(R_{39,*}\) & \( K[x,y,z,w] / (x^2, y^2, z^2 + w^2, xy, xz, xw, yz, yw, zw) \) & A\\
          & 4,1 & \(R_{40}\) & \( K[x,y,z,w] / (x^2, y^2 + z^2, y^2 + w^2, xy, xz, xw, yz, yw, zw) \) & A\\
          & 4,1 & \(R_{41}\) & \( K[x,y,z,w] / (x^2, y^2, z^2, w^2, xy - zw, xz, xw, yz, yw) \) &A\\
          & *4,1 & \(R_{41,*}\) & \( K[x,y,z,w] / ( x^2 + y^2, x^2 + z^2, x^2+ w^2, xy, xz, xw, yz, yw, zw) \) & A\\
          & 5 & \(R_{42}\) & \( K[x, y, z, w, v] / (x, y, z, w, v)^2 \) & A\\
        \hline
    \caption{Local algebras over \( K \) of rank \( \leq 6 \).}
    \label{poonentable}
\end{longtable}

\section{Classification of the Ideal Chomp Game on $\bar{K}$-algebras up to rank $6$}
\label{ClassificationIdealChomp}

We completely classify the Ideal Chomp Game on $\bar{K}$-algebras up to rank 6 and establish a classification of the Ideal Chomp Game on a large portion of Noetherian $\bar{K}$-algebras. Recall Theorem \ref{ClassificationThm1}:

\begin{Thm*}
    [Classification of the Ideal Chomp Game on $\bar{K}$-algebras up to rank 6]
    \label{ClassificationThm2}
    Player A wins on all $\bar{K}$-algebras up to rank 6, except for local $\bar{K}$-algebras isomorphic to 
    \begin{align*}
        R_1 & = \bar{K},\\
        R_4 & = \bar{K}[x,y]/(x,y)^2,\\ R_{12} & =  \bar{K}[x,y]/(xy,x^3,y^3),\\
        R_{13} & = \bar{K}[x,y]/(x^2,xy^2, y^3),\\ R_{17} & = \bar{K}[x,y,z,w]/(x,y,z,w)^2.
    \end{align*}
\end{Thm*}

The proof strategy of Theorem \ref{ClassificationThm2} relies on the Structure Theorem for Artinian Rings (see Theorem 8.7 in \cite{atiyah-macdonald}) which allows us to split any $\bar{K}$-algebra of rank up to 6 into a Cartesian product of local $\bar{K}$-algebras, Poonen's classification of local $\bar{K}$-algebras up to rank 6 (see Table \ref{poonentable}), proving that player B has a winning strategy on $R_1,R_4,R_{12}, R_{13}$ and $R_{17}$ and various reduction steps which show that player A can reduce the game on any $\bar{K}$-algebra of rank up to 6 apart from $R_1, R_4, R_{12}, R_{13}, R_{17}$ to one of these five rings. Thereby, player A is the second player to move on a ring where the second player has a winning strategy. Hence, player A has a winning strategy.

\subsection{Proof of the Classification}

We proceed as follows: In Proposition \ref{Bwins} we prove that player B has a winning strategy on the rings $R_1, R_4, R_{12}, R_{13}$ and $R_{17}$ where Table \ref{poonentable} claims that player B has a winning strategy. In Proposition \ref{AwinsLocal} we then use Table \ref{reducingAlgebras} together with Proposition \ref{Bwins} to also establish all the claimed winning strategies of player A. Altogether these results show Theorem \ref{ClassificationThm2} which gives a complete classification of all Ideal Chomp Games on $K$-algebras up to rank $6$. Some technical statements are deferred from the proofs of these propositions into lemmata.

\begin{Lem}
\label{localalgebras}
    In the game on a local $\bar{K}$-algebra $\bar{K}[x_1,x_2,\dots,x_n]/I$ with the unique maximal ideal $(x_1,\dots,x_n)$ corresponding to the origin $(0,0,...,0)$, any player who plays a polynomial with non-zero constant coefficient loses immediately. In particular, this holds for all the $K$-algebras in Poonen's table \ref{poonentable}.
\end{Lem}

\begin{myproof}
    Let $f$ be a polynomial with non-zero constant coefficient. Then $(0,\dots,0) \notin V(f)$, so in particular $(0,\dots, 0) \notin V((f) + I) = V(f) \cap V(I)$. If $I + (f)$ wasn't the entire ring $\bar{K}[x_1,\dots,x_n]$, then $I + (f)$ would be contained in some maximal ideal $\Fm$. However, $I$ is only contained in $(x_1,\dots,x_n)$ by our assumption that $(x_1,\dots,x_n)$ is the unique maximal ideal in $\bar{K}[x_1,\dots,x_n]/I$. We conclude that $I + (f) = \bar{K}[x_1,\dots,x_n]$, so any player who plays a polynomial with non-zero constant coefficient loses immediately.
\end{myproof}

\begin{Prop}
\label{Bwins}
    Player B has a winning strategy in the Ideal Chomp Game on $R_1, R_4, R_{12}, R_{13}$ and $R_{17}$.
\end{Prop}

\begin{myproof}
    \begin{enumerate}[itemsep = 1ex]
        \item On $R_1$ player B has a winning strategy by Proposition \ref{fieldwin}.

        \item Player B has a winning strategy on $R_4$:\vspace*{1ex}
        
        Suppose player $A$ plays the polynomial $f(x,y) \in R_4$ in the first move. Then $f$ has a representative $\tilde{f} \in \bar{K}[x,y]$ of the form $\tilde{f} = ax + by + c$. By Lemma \ref{localalgebras} we can assume that $c = 0$ as player A would else lose immediately. As $f \neq 0$, we have that $a \neq 0$ or $b \neq 0$. Since the situation is symmetric, let us assume, without loss of generality, that $a \neq 0$. Then player B can play $y$ in his next turn to return the ring $$R_4/(f) + (y) = R_4/(ax + by,y) = R_4/(x,y) \cong \bar{K}$$ to player A who immediately loses in his next move by Proposition \ref{fieldwin}.

        \item Player B has a winning strategy on $R_{12}$:\vspace*{1ex}
        
        Suppose that player A plays $f \in R_{12}\setminus \{0\} = \bar{K}[x,y]/(xy,x^3,y^3)\setminus \{0\}$. Then $f$ has a representative $\tilde{f} \in \bar{K}[x,y]$ of the form $\tilde{f} = a_0 + a_1 x + a_2 y + a_3 x^2 + a_4 y^2$. By Lemma \ref{localalgebras} player A loses immediately if the constant coefficient $a_0$ is non-zero. Hence, we can assume that $\tilde{f}$ has constant coefficient $a_0 = 0$, i.e. $\tilde{f} = a_1x + a_2y + a_3x^2 + a_4y^2$.\vspace*{1ex}
    
        \textit{Case 1:} Suppose that at least one of $a_1,a_2$ is non-zero. Without loss of generality (the situation is symmetric), assume that $a_1 \neq 0$. Since $$\tilde{f} \cdot x = a_1x^2 + a_2xy + a_3x^3 + a_4xy^2 \equiv a_1x^2,$$ we find that $$x^2 \in (xy, x^3,y^3) + (f) = I_0 + (f) = I_1.$$ Hence, $I_1 = (x^2,xy,y^3) + (a_1x + a_2y + a_4 y^2)$ with $a_1 \neq 0$. As a result, it is a valid move for player $B$ to play $y$: the minimal degree of an element of $(x^2, xy, y^3)$ is 2, hence, y would have to lie in $(a_1x + a_2y + a_4y^2)$ which it doesn't since $a_1 \neq 0$. Once player B has played $y$, the game's ideal is $I_2 = (x,y)$, i.e. player A loses in the next turn by Proposition \ref{fieldwin}.\vspace*{1ex}

        \textit{Case 2:} Suppose that $a_1,a_2 = 0$. Then $\tilde{f} = a_3x^2 + a_4y^2$ with at least one of the coefficients $a_3,a_4$ being nonzero. Without loss of generality (the situation is symmetric) assume that $a_3 \neq 0$. Then player B can play $y^2$ in the next turn which reduces the game to the game on $R_4 = \bar{K}[x,y]/(x,y)^2$. This game is won by the second player, i.e. player B, by our analysis of winning strategies on $R_4$.

        \item Player B has a winning strategy on $R_{13}$:\vspace*{1ex}

        Suppose that player A plays $f \in \bar{K}[x,y]/(x^2, xy^2, y^3) \setminus \{0\}$. Then $f$ has a representative of the form $\tilde{f} = a_0 + a_1x + a_2y + a_3 xy + a_4 y^2$. By Lemma \ref{localalgebras} we know that player A loses if the constant term is non-zero, so we may assume $a_0 = 0$. As a preparation we compute:
    
        $$\tilde{f}(x,y) \cdot x = a_2 xy + a_3 x^2y$$
        $$\tilde{f}(x,y) \cdot y = a_1 xy + a_2y^2.$$
    
        \textit{Case 1:} Suppose $a_1 = a_2 = 0$. Then $\tilde{f} = a_3xy + a_4y^2$.\vspace*{1ex}
        
        \hspace*{3mm}\textit{Case 1.1:} If $a_3 \neq 0$, then $y^2 \notin I_1 = I_0 + (\tilde{f})$, so player B may add $y^2$ to the ideal to obtain $I_2 = I_0 + (f) + (y^2) = (x,y)^2$. Hence, the game is reduced to the game on $R_4 = \bar{K}[x,y]/(x,y)^2$. In this case, player B wins by our previous analysis of winning strategies on $R_4$.\vspace*{1ex}
        
        \hspace*{3mm}\textit{Case 1.2:} Otherwise, if $a_3 = 0$, then $\tilde{f} = a_4y^2$. Hence, player B can play $xy$ to reduce the game to the game on $\bar{K}[x,y]/I_0 + (\tilde{f}) + (xy) = \bar{K}[x,y]/(x,y)^2$. Again, player B wins by our previous analysis of winning strategies on $R_4$.\vspace*{1ex}

        \textit{Case 2:} Suppose $a_1 = 0, a_2 \neq 0$. Then by the preparatory computations, $y^2 \in I_1$ and player B can play $x$ since $a_2 \neq 0$. Thereby, player B reduces the game to $$R_{13}/(\tilde{f},x) = \bar{K}[x,y]/(x,y^3,a_2y) = K[x,y]/(x,y) \cong \bar{K}$$\vspace*{1ex} and wins by Proposition \ref{fieldwin}.\vspace*{1ex}

        \textit{Case 3:} Suppose $a_1 \neq 0, a_2 = 0$. Then by the preparatory computations, $xy \in I_1$, but $y^2 \notin I_1$. If player B now plays $y^2$, he reduces the game to $\bar{K}[x,y]/(x,y)^2$ which he wins by our previous analysis of winning strategies on $R_4$.\vspace*{1ex}

        \textit{Case 4:} Finally, suppose $a_1,a_2 \neq 0$.

        Then player B can play $y$. Thereby, he reduces the game to $$R_{13}/I_1 + (y) = \bar{K}[x,y]/(x,y) \cong \bar{K}$$ and wins by Proposition \ref{fieldwin}.

        \item Player B has a winning strategy on $R_{17}$:\vspace*{1ex}

        Any element $f \in R_{17}\setminus\{0\}$ player A plays can be represented by a polynomial of the form $\tilde{f} = a_0 + a_1x + a_2 y + a_3 z + a_4 w$. By Lemma \ref{localalgebras} $a_0 = 0$ if player A tries to avoid losing immediately. Therefore, there exists $i \in \{1,2,3,4\}$ with $a_i \neq 0$. But then taking the quotient $\bar{K}[x,y,z,w]/(x,y,z,w)^2 + (f)$ reduces the rank of $\bar{K}[x,y,z,w]/(x,y,z,w)^2$ over $\bar{K}$ by exactly 1. Looking at Table~\ref{reducingAlgebras} we notice that in all games on $\bar{K}$-algebras with $\rank(R) = 3$, one move suffices to reduce them to the algebras $R_1$ or $R_4$ which are second player wins by what we proved above. Hence, all $\bar{K}$-algebras of rank three are first player wins, making $R_{17} = K[x,y,z,w]/(x,y,z,w)^2$ a second player win.\vspace{1ex}  

        Notice that our proof for the winning strategy on $R_{17}$ only works once we have established that all $\bar{K}$-algebras of rank 3 can be reduced to $R_1$ are $R_4$ by a single move which will happen in Proposition \ref{AwinsLocal} without using the winning strategy on $R_{17}$.
        
    \end{enumerate}
\end{myproof}

For proving that player A has a winning strategy on all the other local $\bar{K}$-algebras up to rank $6$ we need the following auxiliary result:

\begin{Lem}
\label{ringIso1}
     The following rings are isomorphic over any field $K$: $$K[y,z]/(y,z)^2 \cong K[x,y,z]/(x,y,z)^2 +(x + y).$$
\end{Lem}

\begin{myproof}
    Let\vspace*{-5mm} $$\phi: K[y,z]/(y,z)^2 \to K[x,y,z]/(x,y,z)^2 + (x+y),\quad f(y,z) \mapsto f(y,z),$$
    $$\psi: K[x,y,z]/(x,y,z)^2 + (x+y) \to K[y,z]/(y,z)^2,\quad f(x,y,z) \mapsto f(-y,y,z).$$ It is straightforward to check that $\phi$ and $\psi$ are both $K$-algebra homomorphisms. We show that they are two-sided inverses of one another.\\
    Indeed, let $f \in K[x,y,z]/(x,y,z)^2 + (x+y)$. Then $$\phi \circ \psi(f(x,y,z)) = \phi(f(-y,y,z)) = \phi(g(y,z)) = g(y,z) = f(-y,y,z) = f(x,y,z)$$ since $x + y \equiv 0$ in $K[x,y,z]/(x,y,z)^2 + (x+y)$. Conversely, let $f \in K[y,z]/(y,z)^2$. Then $$\psi \circ \phi(f(y,z)) = \psi(f(y,z)) = f(y,z).$$ This concludes the Lemma.
\end{myproof}

The following table summarizes the reduction moves we need for the next Proposition \ref{AwinsLocal}.

\renewcommand{\arraystretch}{1.2}
\setlength\LTleft{0pt plus 1fill}
\setlength\LTright{0pt plus 1fill}
\enlargethispage{3mm}
\vspace*{-3mm}

\begin{table}[h!]
\begin{minipage}[t]{0.48\linewidth}
  \smallskip
    \begin{longtable}{|c|c|c|c|}
        \hline
        \(n\) & ring & move & resulting ring\\
        \hline
        2 & \(R_2\) & \(x\) & \(R_1\)\\
        \hline
        3 & \(R_3\) & \(x\) & \(R_1\)\\
        \hline
        4 & \(R_5\) & \(x\) & \(R_1\)\\
        & \(R_6\) & \(y^2\) & \(R_4\)\\
          & \(R_7\) & \(xy\)& \(R_4\)\\
          & \(R_{7,*}\) & \(x^2\) & \(R_4\)\\
          & \(R_8\) & \(z\) & \(R_4\)\\
        \hline
        5 & \(R_9\) & \(x\) & \(R_1\)\\
          & \(R_{10}\) & \(y^2\) & \(R_4\)\\
          & \(R_{11}\) & \(y^2\) & \(R_4\)\\
          & \(R_{14}\) & \(z\) & \(R_4\)\\
          & \(R_{15}\) & \(z\) & \(R_4\)\\
          & \(R_{15,*}\) & \(z\) & \(R_4\)\\
          & \(R_{16}\) & \(x\) & \(R_4\)\\
        \hline
        6 & \(R_{18}\) & \(x\) & \(R_1\)\\
          & \(R_{19}\) & \(y^2\) & \(R_4\)\\
          & \(R_{20}\) & \(y^2\) & \(R_4\)\\
          & \(R_{21}\) & \(y^3\) & \(R_{12}\)\\
          & \(R_{22}\) & \(x^3\) & \(R_{12}\)\\
          & \(R_{23}\) & \(y^3\) & \(R_{13}\)\\
          & \(R_{24}\) & \(y^3\) & \(R_{13}\)\\
          & \(R_{25}\) & \(xy^2\) & \(R_{13}\)\\
          & \(R_{25,*}\) & \(x^2\) & \(R_{13}\)\\
          & \(R_{25,**}\) & \(y^3\) & \(R_{13}\)\\
          \hline
    \end{longtable}
    \label{reductionMoves}
    \end{minipage}\hfill
    \begin{minipage}[t]{0.48\linewidth}
    \smallskip
    \begin{longtable}{|c|c|c|c|}
        \hline
        \(n\) & ring & move & resulting ring\\
        \hline
          6 & \(R_{26}\) & \(xy\) & \(R_{12}\)\\
          & \(R_{27}\) & \(z\) & \(R_{4}\)\\
          & \(R_{28}\) & \(z\) & \(R_{4}\)\\
          & \(R_{29}\) & \(z\) & \(R_{4}\)\\
          & \(R_{29,*}\) & \(z\) & \(R_{4}\)\\
          & \(R_{30}\) & \(x\) & \(R_4\)\\
          & \(R_{31}\) & \(x\) & \(R_4\)\\
          & \(R_{31,*}\) & \(x\) & \(R_4\)\\
          & \(R_{32}\) & \(z\) & \(R_4\)\\
          & \(R_{33}\) & \(x\) & \(R_{12}\)\\
          & \(R_{34}\) & \(x+y\)& \(R_4\)\\
          & \(R_{34,*}\) & \(x+y\) & \(R_4\)\\
          & \(R_{35}\) & \(x+y\) & \(R_4\)\\
          & \(R_{36}\) & \(z\) & \(R_4\)\\
          & \(R_{36,*}\) & \(z\) & \(R_4\)\\   
          & \(R_{37}\) & \(z\) & \(R_4\)\\
          & \(R_{37,*}\) & \(z\) & \(R_4\)\\
          & \(R_{38}\) & \(w^2\) & \(R_{17}\)\\
          &\(R_{39}\) & \(zw\) & \(R_{17}\)\\
          & \(R_{39,*}\) & \(z^2\) & \(R_{17}\)\\
          & \(R_{40}\) & \(y^2\) & \(R_{17}\)\\
          & \(R_{41}\) & \(xy\) & \(R_{17}\)\\
          & \(R_{41,*}\) & \(x^2\) & \(R_{17}\)\\
          & \(R_{42}\) & \(v\) & \(R_{17}\)\\
        \hline
\end{longtable}
\end{minipage}
\setcounter{table}{1}
\caption{Winning moves for player A on local $\bar{K}$-algebras of rank $\leq 6$.}
\label{reducingAlgebras}
\end{table}

\begin{Prop}
\label{AwinsLocal}
    Player A has a winning strategy on all other rings from Table~\ref{poonentable}, that is on all rings $R_i$ with $i \notin \{1,4,12,13,17\}$.
\end{Prop}

\begin{myproof}
    Consider Table \ref{reducingAlgebras}. It shows that every such $\bar{K}$-algebra $R$ can be reduced to a $\bar{K}$-algebra $R'$ on which the second player has a winning strategy by Proposition \ref{Bwins}. Hence, player A can play the given move in column ``move'' to reduce the current game to a game where the second player, now player A, has a winning strategy. We conclude that player A has a winning strategy on all rings $R_i$ with $i \notin \{1,4,12,13,17\}$ given in Table \ref{poonentable}.

    The only cases where the isomorphism class of the resulting $\bar{K}$-algebra after the move is not immediately clear are $R_{34}, R_{34,*}$ and $R_{35}$. We discuss them here:

    \begin{enumerate}
        \item On $R_{34} = K[x,y,z]/(xy, xz, y^2, z^2, x^3)$ player A has a winning strategy starting with $x + y$. Since
        $$(x+y) \cdot x = x^2 + xy \equiv x^2 \text{ and } (x+y) \cdot z \equiv yz,$$ 
        we have $x^2, yz \in I_0 + (x+y).$ Hence, $$R_{34}/(x+y) = \bar{K}[x,y,z]/(x,y,z)^2 + (x+y),$$ so by Lemma \ref{ringIso1} we conclude that player A reduces the game to a game isomorphic to $K[y,z]/(y,z)^2 \cong R_4$ where the second player (then player A) has a winning strategy.

        \item On $R_{34,*} = \bar{K}[x,y,z]/(xy, xz, yz, y^2-z^2, x^3)$ player A has a winning strategy starting with $x + y$ and, thereby, reducing the game to
         \begin{align*}
             R_{34,*}/(x+y) = \bar{K}[x,y,z] / (xy, xz, yz, y^2 - z^2, x^3) + (x+y) =\\= \bar{K}[x,y,z]/(x,y,z)^2 + (x+y) \cong \bar{K}[y,z]/(y,z)^2 \cong R_4
         \end{align*} by Lemma \ref{ringIso1} and since in $R_{34,*}$ we have $$(x+y) \cdot y = xy + y^2 = y^2 \text{ and } (x+y)x = x^2 + xy = x^2.$$

         \item On $R_{35} = \bar{K}[x,y,z]/(xy, xz, yz, x^2 + y^2 - z^2)$ player A has a winning strategy starting with $x+y$. Since $$(x+y)x = xy + x^2 \equiv x^2 \text{ and }(x+y)y = xy + y^2 \equiv y^2,$$ we have that $x^2, y^2, z^2 \in I_0 + (x+y)$. Hence, by Lemma \ref{ringIso1} $R_{35}/(x+y) = \bar{K}[x,y,z]/(x,y,z)^2 + (x+y) \cong \bar{K}[y,z]/(y,z)^2 \cong R_4$.
    \end{enumerate}

     We also verified all the reductions given in Table \ref{reducingAlgebras} using Sage. The Sage Worksheet can be found under the link provided in the bibliography \cite{karl}.
\end{myproof}

\begin{Cor}
    Table \ref{poonentable} classifies the Ideal Chomp Game on all local $\bar{K}$-algebras up to rank $6$.
\end{Cor}

We conclude this section by classifying the Ideal Chomp Game on all (not necessarily local) $\bar{K}$-algebras up to rank $6$ and provide partial results on $\bar{K}$-algebras of Krull dimension~0 and higher rank.

\begin{Prop}
    Player A has a winning strategy on every non-local $\bar{K}$-algebra of finite rank.
\end{Prop}

\begin{myproof}
    Since finite-rank $\bar{K}$-algebras are Artinian rings, the Structure Theorem for Artinian Rings (see e.g. Theorem 8.7 in \cite{atiyah-macdonald}) implies that any finite-rank $\bar{K}$-algebra can be written as a finite direct product of local finite-rank $\bar{K}$-algebras. Since the $\bar{K}$-algebras are non-local, this direct product is non-trivial. Hence, player A has a winning strategy on all non-local $\bar{K}$-algebras of finite rank by Proposition \ref{CartesianProdCor}.
\end{myproof}

Returning to our main classification result, we obtain:

\begin{Thm}
    [Classification of the Ideal Chomp Game on $\bar{K}$-algebras up to rank 6]
    \label{ClassificationThm}
    Player A wins on all $\bar{K}$-algebras up to rank 6, except for local $\bar{K}$-algebras isomorphic to
    \begin{myitemize}
        \item $R_1 = \bar{K}$,
        \item $R_4 = \bar{K}[x,y]/(x,y)^2$, \item $R_{12} =  \bar{K}[x,y]/(xy,x^3,y^3)$,
        \item $R_{13} = \bar{K}[x,y]/(x^2,xy^2, y^3)$ and \item $R_{17} = \bar{K}[x,y,z,w]/(x,y,z,w)^2$.
    \end{myitemize}
\end{Thm}

\begin{myproof}
    By Corollary \ref{CartesianProdCor} player A wins on all finite-rank $\bar{K}$-algebras that are not local themselves, including those up to rank $6$. By Proposition \ref{Bwins} and Proposition \ref{AwinsLocal} $\bar{K}$-algebras isomorphic to $R_1, R_4, R_{12}, R_{13}$ and $R_{17}$ are the only local $\bar{K}$-algebras up to rank 6 where player B has a winning strategy. This concludes the proof.
\end{myproof}

\section{Extensions and Outlook}

After our complete classification of $\bar{K}$-algebras of Krull dimension 0 up to rank~$6$ and our partial results on $\bar{K}$-algebras of Krull dimension 0 and higher rank in Section \ref{ClassificationIdealChomp}, we will provide explicit examples and ideas on how to approach the Ideal Chomp Game on $\bar{K}$-algebras of higher Krull dimensions. Furthermore, we will state questions that are still open to the author's best knowledge along the way.

\subsection{Henson's Theorem}
\label{sec:henson}

In 1970 Henson solved the Ideal Chomp Game on a class of rings that are close to being principal ideal domains, extending the result that player A has a winning strategy on principal ideal domains that are not fields (Proposition \ref{principalmaxidealwin}) to these rings. We discuss this result and give an example of a class of rings on which the Ideal Chomp Game is solved using Henson's result. Furthermore, we discuss its limitations.

\begin{Def}
\label{hensoncondition}
    A ring $R$ satisfies Henson's condition (with respect to $a\in R$) if $R$ is an integral domain and there exists an element $a \in R$ such that $R/(a)$ is a PID, but not a field.
\end{Def}

\begin{Thm}
[Henson's Theorem]
\label{hensontheorem}
    Let $R$ be a ring satisfying Henson's condition. Then the first player has a winning strategy in
    the Ideal Chomp Game on $R$, beginning with the move $a^2$.
\end{Thm}

\begin{skproof}
    The main idea of the proof is to define so-called \emph{special ideals}. These are ideals $I \ideal R/(a^2)$ that can be written as a power of some maximal ideal, i.e. $I = \Fm^b$ for some $b \in \BN_{\geq 1}$. Using some more technical lemmata, one shows that player A can always reach such a special ideal when playing correctly while prohibiting player B from reaching such a special ideal. Since any maximal ideal is in particular special, player A will be the player to build some maximal ideal. Hence, player B loses the game. The complete proof can be found in \cite{henson}.
\end{skproof}

Some examples of rings that satisfy Henson's condition w.r.t. $x$ are $K[x,y]$ for an arbitrary field $K$, $\BZ[x]$ or more generally $R[x]$ for any principal ideal domain $R$.

We now discuss an explicit winning strategy for player A on $\bar{K}[x,y]$. The needed background on primary ideals can be found in \cite{atiyah-macdonald}.

\begin{Ex}[A winning strategy for player A on $\bar{K}\lbrack x,y\rbrack$]
    Player A starts by playing $a_1 = x^2 \in \bar{K}[x,y]$ according to Henson's Theorem \ref{hensontheorem}. Player $B$ then plays $$a_2 = f(x,y) = p(y) + x\cdot q(y) \in \bar{K}[x,y]/ (x^2).$$
    Let $R = \bar{K}[x,y]/(x^2)$, for all $n \geq 2$ the ideal $I_n = (a_2, \dots,a_n) \ideal R$ is always understood as an ideal in $R$ and all the $a_i$'s satisfy $a_i \in R \setminus (a_2,\dots, a_{i-1})$.\\
    Notice that all maximal ideals of $R$ are of the form $\Fm_b = (x,y-b)$ for some $b \in \bar{K}$. We call an ideal \emph{special} if it is of the form $\Fm_b^k$ for some $b \in \bar{K}$ and some $k \geq 1$. The special ideals in $R$ are precisely of the form $$\Fm_b^k = 
        ((y-b)^{k}, x(y-b)^{k-1}).$$

    In particular, $I_2 = (a_2) = (f)$ is not special as it is principal while all special ideals are not.

    Let $$f(x,y) = p(y) + x \cdot q(y) = c_p \cdot \prod_{i=1}^n (y - b_i)^{s_i} + x \cdot c_q \cdot \prod_{i = 1}^n (y - b_i)^{t_i}$$ where $c_p, c_q \in \bar{K}$. We discuss a special case where some primary ideal over $f$ is easy to find. The general approach to finding a response of player A to a polynomial player B plays can be found in \cite{henson}.

    Suppose there exists $b_i$ with $y-b_i \mid p,q$. Let $r = \min\{s_i,t_i\}$. Then the $\Fm_{b_i}$-primary ideal over $f$ is $\left((y-b_i)^{r}\right)$ and player A can choose $$a_3 = (y-b_i)^r + x(y-b_i)^{r-1}$$ to build the ideal $$I_3 = (f,a_3) = \Fm_{b_i}^r.$$ Hence, player A was able to build a special ideal. Player A continues to do so in the following rounds and finally wins the game.

\end{Ex}

Let us discuss one more explicit example:

\begin{Ex}
    
    Let player A start with $a_1 = x^2$ and player B play $a_2 = x(y-1) + (y-2)$. Then player A will respond with $a_3 = y-2$ which yields the ideal $$I_3 = (x(y-1), (y-2), x^2) = (x, y-2),$$ since $x = x(y-1) - x(y-2) \in I_3$. As this is a maximal ideal, player A wins.
\end{Ex}

Notice, however, that the use of Henson's theorem is limited when exploring the Ideal Chomp Game on higher dimensional rings:

\begin{Lem}
    Let $R$ be a Noetherian local ring of Krull dimension larger than or equal to three. Then $R$ does not satisfy Henson's condition.
\end{Lem}

\begin{myproof}
    By corollary 11.18 from \cite{atiyah-macdonald}, the dimension of $R/(x)$ for $x$ in the unique maximal ideal $\Fm \ideal R$ is $$\dim(R/(x)) = \dim(R) - 1.$$ Since principal ideal domains are of Krull dimension less than or equal to $1$, we conclude.
\end{myproof}

\subsection{The Ideal Chomp Game on higher dimensional $\bar{K}$-Algebras}
\label{weiteresKalgebren}

We already classified a large portion of Ideal Chomp Games on $\bar{K}$-algebras of Krull dimension 0 in Theorem \ref{ClassificationThm}. What remains towards a complete classification is a classification of the Ideal Chomp Game on local $\bar{K}$-algebras of finite rank $\geq 7$ and of local $\bar{K}$-algebras of Krull dimension~$\geq 1$. In the following, we discuss some partial results in the Ideal Chomp Game on $\bar{K}$-algebras with Krull dimension~$\geq 1$.

\subsubsection{$\bar{K}$-Algebras of Krull Dimension 1}

Besides Henson's Theorem (see Theorem \ref{hensontheorem}), there is another result on winning strategies in the Ideal Chomp Game on $\bar{K}$-algebras of Krull dimension 1:

\begin{Prop}
\label{ellipticCurves}
    Let $f = 0$ be an affine Weierstrass equation in $\bar{K}[x,y]$. Then player B has a winning strategy on $\bar{K}[x,y]/(f)$.
\end{Prop}

\begin{myproof}
    See Brandenburg's Paper ``Algebraic games -- Playing with groups and rings'' (\cite{brandenburg}), Proposition 5.7.
\end{myproof}

For $\bar{K}$-algebras of Krull dimension 1, a classification of the Ideal Chomp Game on them is missing.

\begin{Ques}
    Can we classify the Ideal Chomp Game on $\bar{K}$-algebras of Krull dimension 1?
\end{Ques}

\subsubsection{$\bar{K}$-Algebras of Krull Dimension 2}

For the Ideal Chomp Game on $\bar{K}$-algebras of Krull Dimension 2 we don't have any results that classify the game on a larger portion of these algebras. Hence, we will just look at one specific example.

\begin{Ex}
    In the Ideal Chomp Game on $\bar{K}[x,y]$ player A has a winning strategy. This can be obtained from both Henson's Theorem (see Theorem \ref{hensontheorem}) and the proposition on winning strategies on elliptic curves (see Proposition \ref{ellipticCurves}).
\end{Ex}

For most Ideal Chomp Games on $\bar{K}$-algebras, or more generally on most rings, of Krull dimension 2 it is not known which player has a winning strategy. In particular, the following question is open for any choice of your favorite ring that does not satisfy Henson's condition or is a non-trivial Cartesian product of rings, e.g. for $K[x,y,z]/(z^{a})$ for $a \in \BN_{\geq 2}$.

\begin{Ques}
    Who has a winning strategy in the Ideal Chomp Game on your favorite (local) $\bar{K}$-algebra/ring of Krull dimension 2? 
\end{Ques}

\subsubsection{$\bar{K}$-Algebras of Krull Dimension $\geq 3$}
\label{sec:KDgeq3}

By Proposition \ref{CartesianProdCor} player A has a winning strategy on rings of the form $R = \prod_{i=1}^{n} R_i$ with $R_i$ of Krull dimension 1 and $n \geq 2$. These rings are of Krull dimension~n. In particular, there are examples of rings of Krull dimension n where we know which player has a winning strategy. For player A we even have examples of rings for any given Krull dimension $\geq 1$ that cannot be written as a nontrivial Cartesian product of rings and where player A has a winning strategy:

\begin{Ex}
    Given $n \in \BN$, consider the scheme $S$ we receive by gluing together the affine space $\BA^{n-1}$ with $\BA^1$ at the origin in $\BA^n$. Then the coordinate ring $R \subseteq \bar{K}[z_1,z_2,z_3,z_4]$ of $S$ cannot be written as a nontrivial Cartesian product since $S$ is connected. Furthermore, player A has the winning strategy to play $z_4 - 1$ in his first move to reduce $R$ to $R/(z_4-1) \cong K$ where the isomorphism follows from the fact that $S \cap \{z_4 = 1\} = \{(0,0,0,1)\}$ is a single point.
\end{Ex}

Since there are so many ways for the first player to reduce $\bar{K}[x_1,\dots,x_n]$ to a $\bar{K}$-algebra of Krull dimension $n-1$, we conjecture the following regarding this open question:
\begin{Conj}
\label{IdealChompGameConjecture}
    Player A has a winning strategy on $\bar{K}[x_1,\dots,x_n]$ for any $n \geq 1$. 
\end{Conj}

Similarly, we ask the question whether or not a similar statement holds for player B. We expect this question to be more difficult than the previous one.

\begin{Ques}
    Can we find a ring $R_n$ in any Krull dimension $n \in \BN$ where player~B has a winning strategy?
\end{Ques}

Finally, we show that there are infinitely many Krull dimensions in which player A and player B, respectively, have a winning strategy.

\begin{Lem}
    There exists a ring $R_n$ in every Krull dimension $n \in \BN_{\geq 0}$ such that player A has a winning strategy in the Ideal Chomp Game on $R_n$.
\end{Lem}

\begin{myproof}    
    Since the ring $$R_n = \bar{K}[x_0]/(x_0^2) \times  \bar{K}[x_1,\dots,x_n]$$ is a non-trivial cartesian product of Notherian rings, player A has a winning strategy on each $R_n$ by Proposition \ref{CartesianProdCor}. As $\bar{K}[x_0]/(x_0^2)$ has Krull dimension 0 and $\bar{K}[x_1,\dots,x_n]$ has Krull dimension $n$ by ..., $R_n$ has Krull dimension $n$. This concludes the proof.
\end{myproof}

\begin{Lem}
    There exist infinitely many Krull dimensions $n \in \BN$ such that there exists a ring $R_n$ of Krull dimension $n$ where player B has a winning strategy.
\end{Lem}

\begin{myproof}
    Assume there was a Krull dimension $n \in \BN$ such that player A has a winning strategy on every ring of Krull dimension $n$. Then player B has a winning strategy in the Ideal Chomp Game on $\bar{K}[x_1,\dots, x_{n+1}]$ since any element player A can choose in their first turn reduces the ring to a ring of Krull dimension $n$. Hence, player B has a winning strategy on infinitely many rings of pairwise different Krull dimension.
\end{myproof}

\pagebreak

%%%%%%%%%%%%%%%%%%%%%%%%%%%%%%%%%%%%%%%%%%%%%%%%%%%%%%%%%%%%%
% Bibiliography, Appendix and Eigenständigkeitserklärung

\section*{Acknowledgements}

This paper contains the results I found in my Master's thesis and is mainly an abbreviation thereof. I would like to sincerely thank my supervisor Johannes Schmitt for his endless support both while writing my Master's thesis and in the publishing process.

You can find my Master's thesis in the ETH Research Collection under the following link:
\url{https://www.research-collection.ethz.ch/entities/publication/f32e7185-384b-42fe-bb3a-ff749d48e0ee}.

\vfill

Leopold Karl, Institut für Diskrete Mathematik \& Geometrie, FB 1 Algebra, Technische Universität Wien, Wiedner Hauptstraße 8-10/104, 1040 Vienna, Austria.

E-mail address: leopold.karl@tuwien.ac.at


\begin{thebibliography}{Vin}
\addcontentsline{toc}{section}{References}

\bibitem[AMD]{atiyah-macdonald} M. F. Atiyah, I. G. MacDonald, Introduction to Commutative Algebra, CRC Press, Boca Raton, 2018.

\bibitem[Bra]{brandenburg} M. Brandenburg, Algebraic games -- Playing with groups and rings, 2017, \url{https://arxiv.org/abs/1205.2884}.

\bibitem[Gal]{gale} D. Gale, A Curious Nim-Type Game, The American Mathematical Monthly, Vol. 81, No. 8 (Oct., 1974), pp. 876--879.

\bibitem[Gat]{gathmann} A. Gathmann, Commutative Algebra, Class Notes TU Kaiserslautern 2013/14, \url{https://agag-gathmann.math.rptu.de/en/commalg.php}.

\bibitem[Hen]{henson} C.W. Henson, Winning strategies for the Ideal Game, The American Mathematical Monthly, Vol. 77, No. 8 (Oct., 1970), pp. 836--840.

\bibitem[Karl]{karl} L. Karl, SageWorksheet: Computational verification of reduction moves in the Ideal Chomp Game, \url{https://leopoldkarl.github.io/SageWS.pdf}.

\bibitem[Poo]{poonen} B. Poonen, Isomorphism types of commutative algebras of finite rank over an algebraically closed field, Contemporary Mathematics, Vol. 463, 2008, pp. 111--120.

\bibitem[Sch]{schuh} F. Schuh: Spel van delers. Nieuw Tijdschrift voor Wiskunde 39 (1952), pp.~299--304.

\bibitem[Zer]{zermelo} E. Zermelo, Über eine Anwendung der Mengenlehre auf die Theorie des Schachspiels, 1913, \url{https://webdocs.cs.ualberta.ca/~hayward/396/asn/zermelo.pdf}.

\end{thebibliography}
\end{document}